\input amstex
\documentstyle {amsppt}

\magnification1000

\centerline{\bf On the multiplicities of families of complex
hypersurface-germs}
\centerline{\bf  with constant Milnor number}

\bigskip
\centerline{\bf  Camille Pl\'enat and David Trotman}

\medskip
\centerline{\bf  April 17th, 2012}
\bigskip

{\bf Abstract.}
{\it We show that the possible drop in multiplicity in an analytic family
$F(z,t)$ of complex analytic hypersurface singularities with constant
Milnor number is controlled by the powers of $t$. We prove 
equimultiplicity of  $\mu$-constant families of the form $f + tg + t^2h$
if the singular set of the tangent cone of $\{ f = 0 \}$ is not contained
in the tangent cone of $\{ h = 0 \}$. }

\bigskip

{\bf 1. Background.}

\medskip

Let  $F : ({\Bbb C}^n \times {\Bbb C},0\times {\Bbb C}) \longrightarrow
({\Bbb C},0)$ be a holomorphic function in $z_1, \dots, z_n$ and $t$.

We study the following conjecture, stated implicitly by Teissier in 1974
in his Arcata survey [18] as well as at the beginning of his Carg\`ese
paper [17], and which implies a parametrised version of Zariski's problem
[24] about the topological invariance of the multiplicity (Conjecture 1.2
below).

\bigskip

{\bf Conjecture 1.1 (Teissier 1972 [17]).} {\it If  for  $t$ in some
neighbourhood $U$ of $0$ in $ {\Bbb C}$,  each function $F(.,t)$ has an
isolated singularity at the origin with the same Milnor number $\mu$, then
the functions $F(.,t)$ have the same multiplicity at $0$ for $t \in U$.}

\medskip
Teissier made the stronger conjecture at Carg\`ese in 1972 that
$\mu$-constancy implies that the Whitney conditions hold for $(F^{-1}(0),
0 \times {\Bbb C})$. (The same conjecture was made by L\^e D\~ung Tr\'ang
and Ramanujam in [11], published in 1976, although submitted in June
1973.) This turned out to be false as first illustrated by the famous
examples of Brian\c con and Speder [5].   Thus Conjecture 1.1 may also be
considered as a conjecture of Teissier which remains open. It has two
corollaries, as follows.

\bigskip
{\bf Conjecture 1.2 (Zariski's problem for families).} {\it  Families of
complex analytic hypersurfaces with isolated singularities of constant
topological type are equimultiple.}

\medskip
{\bf Proof.} This would follow from Conjecture 1.1 because the Milnor
number is a topological invariant (Milnor [13], Teissier [17]). \qed

\medskip
{\bf Conjecture 1.3} {\it Bekka $(C)$-regular families of complex
hypersurfaces are equimultiple.}
\medskip
{\bf Proof.} Use the analogue of the Thom-Mather isotopy theorem  for
$(C)$-regularity as proved by Bekka in his thesis [1], together with
Conjecture 1.2. \qed

\medskip
Equimultiplicity was established in the case of Whitney regularity for
general complex analytic varieties by Hironaka [9], and with a different
proof,  for Whitney regularity of  families of complex analytic
hypersurfaces  by Brian\c con and Speder [6].  The proof of Brian\c con
and Speder was first extended to arbitrary complex analytic varieties by
Navarro Aznar in [14], and the result is a  special case of Teissier's
general characterisation [19] of Whitney regularity in terms of
equimultiplicity of  polar varieties.

\medskip
Conjecture 1.3 for the stronger hypothesis of weak Whitney regularity
(defined by Bekka and Trotman [2], weak Whitney regularity  is weaker than
Whitney regularity but stronger than  $(C)$-regularity [3])  was proved
directly in 2010 [22] by the second author and  Duco van Straten,  i.e.
weak Whitney regularity  implies equimultiplicity for families of complex
hypersurfaces.

\medskip
Conjecture 1.2 is still unproved, as is Conjecture 1.3. It is also unknown
whether constant topological type implies $(C)$-regularity. It was shown
recently by Bekka and Trotman  [4] that $(C)$-regularity is in general
weaker than weak Whitney regularity for the 1975 quasi-homogeneous
examples of Brian\c con and Speder [5]. No example is currently known of a
weakly Whitney regular complex analytic stratification not also satisfying
Whitney regularity, while the equivalence of the two conditions has only
been proved in the classical case of a family of plane curves. For this,
use that weak Whitney regularity implies $(C)$-regularity [3], which
implies topological triviality by [1], and hence constant Milnor number,
and it is well-known that Whitney regularity is equivalent to constancy of
the Milnor number for families of plane curves [18].

\medskip

{\bf The L\^e-Saito-Teissier criterion for $\mu$-constancy.}

\medskip
According to L\^e and Saito [10] and Teissier [17], the constancy of the
Milnor number of $F(.,t)$  is equivalent to $F$ being a Thom map, i.e. 
equivalent to the $(a_F)$ condition being satisfied. This can be
reformulated as saying that

$${\vert F_ t \vert \over \vert \vert F_z \vert \vert} \rightarrow 0 \,\,
\text{as} \,\, (z,t) \rightarrow (0,0)$$

\noindent where $F_t$ is notation for $\partial F / \partial t$, $F_z$ is
notation for $(\partial F / \partial z_1, \dots , \partial F / \partial
z_n)$,  $\vert .  \vert$ denotes the modulus of a complex number and
$\vert \vert . \vert \vert$ denotes the hermitian norm on ${\Bbb C^{n}}$.

In this paper we use this criterion for  constancy of the Milnor number 
to determine some situations when equimultiplicity holds (Propositions 1.1
and 3.2), and to reduce possible jumps in the multiplicity (Propositions
2.1 and 2.2).

\bigskip

Write $F(z,t) = f(z) + \sum\limits_{k\geq1} t^k g_k(z)$.

Then $F_t = \sum\limits_{k\geq 1} k t^{k-1} g_k$, and $F_z = f_z +
\sum\limits_{k \geq 1}t^k (g_k)_z$.

Let $m(h)$ denote the multiplicity of a function $h$ at $0$.

Due to its upper semicontinuity the multiplicity is non-constant if and
only if $m = m(f) >  m_1 = \min\limits_{k\geq 1} m(g_k)$.

\bigskip
{\bf Proposition 1.1} {\it If $F(z,t) = f(z) + tg(z)$ is a $1$-parameter
family of isolated complex analytic hypersurface singularities whose
Milnor numbers are constant, then the multiplicity at $0$ of $g$ is
greater than or equal to the multiplicity at $0$ of $f$.}

\medskip
 {\bf Proof.}

Suppose that $m(f) = m > m_1 = m(g)$.

Consider  analytic arcs  $\gamma : ({\Bbb C}, 0)  \longrightarrow ({\Bbb
C^{n+1}}, 0)$, $\gamma(u) = (z(u), t(u))$, such that $\gamma(0) = 0 \in
{\Bbb C^{n+1}}$.  We must find an arc $\gamma$ such that

 $${\vert F_ t (\gamma(u)) \vert \over \vert \vert F_z(\gamma(u)) \vert
\vert} \not\rightarrow 0 \,\, {\text{as}} \,\, u \rightarrow 0.$$

For any analytic function $Q : {\Bbb C}^n \times {\Bbb C} \longrightarrow
{\Bbb C}$, write  $V(Q)$ for the order in $u$ at $0$ of $Q \circ \gamma :
{\Bbb C} \longrightarrow {\Bbb C}$, and for any analytic function $P :
{\Bbb C}^n \longrightarrow {\Bbb C}$ write $v(P)$ for the order in $u$ at
$0$ of $P \circ \pi_z \circ \gamma$ , where $\pi_z$ is the projection on
the $n$ first coordinates.

We must  choose an analytic arc $\gamma$ such that $V(F_t) - min \{
V(\partial F/ \partial{z_i})\} \leq 0$.

Now $V(F_t) - min \{ V(\partial F/ \partial{z_i})\}
=  v(g) - min \{V(\partial f/ \partial {z_i} + t \partial g/ \partial
{z_i})  \}$.

Let $\gamma(u) = (uz_0,0)$ where $z_0 \in  {\Bbb C}^n - \{ 0 \} $.

Then
$V(F_t ) - min \{ V(\partial F/ \partial{z_i})\}  =  v(g) - min
\{v(\partial f/ \partial {z_i})  \}.$
For $z_0$ sufficiently general, the right-hand side is
$v(g) - (v(f) - 1) = m_1 - m + 1 \leq 0$,
because  $m_1 < m$.

Thus
$V(F_t) - min \{ V(\partial F/ \partial{z_i})\} \leq 0,$
contradicting the hypothesis that $\mu$ be constant, using the
L\^e-Saito-Teissier characterisation.
$\qed$

\medskip

{\bf Remark 1.2.}  The result of Proposition 1.1 was discovered by the
second author during the academic year 1976-77 and announced in a talk
given in March 1977 [20] as one of the weekly A'Campo-MacPherson
singularity seminars at the University of Paris 7. The text of this talk
was included in the second author's Th\`ese d'\' Etat [21] defended at
Orsay in January 1980. The result  was rediscovered by Gert-Martin Greuel
in 1986 [8], and used to prove Teissier's conjecture in the case of
quasi-homogeneous, and semi-quasihomogeneous functions $f$.

\medskip

{\bf Remark 1.3.} Parusinski [16] has proved that a $\mu$-constant family
of the form $f + tg$ has constant topological type by integrating an
appropriate vector field, with an argument which works for all $n$. L\^e
D\~ung Tr\'ang and Ramanujam [11] proved that for a $\mu$-constant family
of complex hypersurfaces defined by $\{ F(z,t) = 0 \}$, the hypersurfaces 
$\{ F(z,t) = 0 \} \cap ( {\Bbb C^n} \times \{ t \} )$ have constant
topological type when $n \neq 3$.

\medskip

We thank T. Gaffney and M. Oka for their interest and comments on this
work, and also thank the anonymous referee for suggestions leading to
improvements in the paper.

\bigskip

{\bf 2. Controlling multiplicity.}

\medskip

{\bf Proposition 2.1.} {\it  If $F(z,t) = f(z) + tg(z) + t^2 h(z)$ is an
analytic $1$-parameter family of isolated hypersurface singularities with
Milnor number $\mu$ constant, then the multiplicity at the origin of $g$
is greater than or equal to the multiplicity $m$ at the origin of $f$, and
the multiplicity at $0$ of $h$ is greater than or equal to $m - 1$.}
\medskip
{\bf Proof.}
 Because

 $${ \vert F_t \vert \over \vert\vert F_z \vert\vert} = {\vert g + 2t h
\vert \over \vert\vert f_z + tg_z + t^2 h_z \vert \vert } ,$$

\noindent  it follows that on a generic curve of the form $(u z_0,0)$
with $z_0 \neq 0$, $V(F_t) - V(F_z) = m(g) - v(f_z) = m(g) - m + 1$.
 Hence if $m(g) - m +1 \leq 0,$ i.e. $m(g) < m$,  we obtain a
contradiction to the hypothesis that the Milnor number remains constant.

 Thus we obtain that the coefficient $g$ of $t$ has multiplicity
 $m(g) \geq m = m(f)$.

 Suppose that  $m(h) \leq m - 2$.

 On a generic curve of the form $(uz_0, ut_0)$, with both $z_0 \neq 0$
and $t_0 \neq 0$, if $\Delta = {\vert F_t \vert \over \vert \vert F_z
\vert \vert}$, then $\Delta \sim {\vert 2t h \vert \over \vert\vert f_z
+ t^2 h_z \vert\vert} $.

 Hence $V(\Delta) = 1 + m(h) - min \{ m - 1, 2 + m(h) - 1\} = 1 + m(h) -
(m(h) + 1) = 0$.

This again contradicts the hypothesis that the Milnor number of the family
$F(.,t)$ is constant, proving  that $m(h)$ is at least  $m - 1$. \qed

\bigskip
Now we generalize to arbitrary deformations of $f$ which are analytic in
$t$.

\bigskip

{\bf Proposition 2.2.} {\it If  the family $F(z,t) = f(z) + tg_1(z)
+t^2g_2(z) + t^3g_3(z) + \dots $ has constant Milnor number at $z=0$ , and
$f$ has multiplicity $m$ at the origin,  then $m(g_1) \geq m$, $m(g_2)
\geq m - 1, \dots,$ and in general $ m(g_r) \geq m - r+1$.}
\medskip

{\bf Proof.}

Here,
$V(F_t) - V(F_z) = \min \limits_{k \geq 1} \{ (k-1) + v(g_k) \} - V( f_z
+ \sum\limits_{k \geq 1} t^k (g_k)_z)$,
assuming that we are on a generic arc for which there is no cancellation
of terms in the expression for $F_t$.

In particular, $V(F_t)= \min \{ m(g_1) , m(g_2)+1 ,\dots, m(g_r)+r-1 \} $,
for some $r$ such that $r \leq m(g_1)$ (because each $m(g_k) > 1$), and

$V(F_z)= \min \{ m(f)-1, m(g_1)-1+1 , m(g_2)-1+2 , \dots,  m(g_r)-1+r \}$

$\hskip11mm = \min \{ m - 1, V(F_t) \}$.

But the family $F(z,t)$ has constant Milnor number, so that $
V(F_t)>V(F_z)$, by the L\^e-Saito-Teissier theorem.
The conclusion follows.
\qed

\medskip

It is interesting to compare the previous result with the following,
proved by Greuel in 1986 [8].
Observe the extra restrictions Greuel imposed on the $\{ g_i \}$.

\medskip
{\bf Proposition 2.3 (Greuel).} {\it Let $\lambda_j : ({\Bbb C},0)
\rightarrow ({\Bbb C},0)$ and $g_j : ({\Bbb C}^n, 0) \rightarrow ({\Bbb
C},0),  j = 1, \dots, r$,  be holomorphic functions such that $\lambda_j
\neq 0$ and the initial forms of $g_j$ are ${\Bbb C}$-linearly
independent. Assume that

$$F(z,t) = f(z) + \Sigma^{r}_{j=1} \lambda_j(t) g_j(z)$$

\noindent is a $\mu$-constant unfolding of $f$. Then $\nu(\lambda_j) +
m(g_j) > m(f)$ for all $j = 1, \dots, {r}$, where $\nu(\lambda_j)$ denotes
the order in $t$ of $\lambda_j(t)$.}

\bigskip
{\bf Remarks 2.4.}

\medskip
1. Let  $F(z,t) = f(z) + t g_1(z) +t  g_2(z)$
be a $\mu$-constant family, where $g_1(z) = h_{m-1}(z) +  k_1(z)$  and
$g_2(z)  = - h_{m-1}(z) + k_2(z)$,  $m = m(f)$, $h_{m-1}$ is a homogeneous
polynomial of degree $m-1$, and each of $k_1$ and $k_2$ has order greater
than equal to $m$. In this form Greuel's Proposition 2.3 does not apply
because the hypothesis of linear independence does not hold. Also the
conclusion fails to hold because $\nu(t) + m(g_i) = 1 + m - 1 = m = m(f)$,
for $i = 1, 2$. So this explains why Greuel imposes the extra condition of
linear independence of the initial forms of the $g_i$'s.
However if we group together the coefficients of $t$, we can express
$F(z,t)$ as $F(z,t) = t (g_1(z) + g_2(z))$, and Proposition 1.1 applies.

\medskip

2.  Here is a less trivial example - in the previous example $\lambda_1 =
\lambda_2$. Now we provide an example where the $\lambda_i$'s are
distinct. Let

$$F(z,t) = f(z) + (t^3 + at^2)g_1(z) + (t^4 + t^2)g_2(z) + (at^4 -
t^3)g_3(z)$$

\noindent be a $\mu$-constant family, where $g_1(z) = h_{m-2}(z) + 
k_1(z)$, $g_2(z)  = - a h_{m-2}(z) + k_2(z)$, $g_3(z)  = h_{m-2}(z) +
k_3(z)$, $m = m(f)$, $h_{m-2}$ is a homogeneous polynomial of degree
$m-2$, and each of $k_1$, $k_2$, $k_3$  has order greater than equal to
$m$. Then again Proposition 2.3 does not apply. Also the conclusion of
Proposition 2.3 fails to hold because
$$\nu(\lambda_1) + m(g_1) = 2 + m - 2 = m,  {\text{and}} \,\, 
\nu(\lambda_2) + m(g_2) = 2 + m - 2 = m.$$

In this example the $\lambda_i$'s are distinct. We can easily construct
similar more complicated examples all illustrating why Greuel's extra
condition is required.

Grouping together the coefficients of powers of $t$ in this example as in
our Proposition 2.2 we find

$$F(z,t) = f(z) + t^2 (a k_1 + k_2) + t^3 (k_1- k_3) + t^4 (k_2 + a k_3).$$

It may be the case that Greuel's Proposition 2.3 still does not apply.
This will occur if the initial terms of $(ak_1 + k_2)$, $(k_1 - k_3)$ and
$(k_2 + a k_3)$ are linearly dependent. Let $k_i = h_r + p_i$ where
$h_r(z)$ is a homogeneous polynomial of degree $r$ and $p_i(z)$ is a
polynomial of order at least $r + 1$, for $i = 1$ and $i = 3$.  Then if
$k_2$ is any homogeneous polynomial of degree $r$ which is not equal to $-
a h_r$, then the coefficients of $t^2$ and $t^4$ will have initial forms
of degree $r$  which are the same, hence linearly dependent, so that the
hypothesis of Proposition 2.3 is not satisfied. Our Proposition 2.2 allows
us to prove that $r \geq m - 1$, so that the conclusion of Proposition 2.3
($2 + r  > m $), does hold.

\medskip

3. Here is the point in the proof of Proposition 2.2 that allows us to
proceed without Greuel's  restriction on the initial forms of the $g_i$.

We say:
$V(F_t) - V(F_z) = \min \limits_{k \geq 1} \{ (k-1) + v(g_k) \} - V( f_z
+ \sum\limits_{k \geq 1} t^k (g_k)_z)$,
assuming that we are on a generic arc for which there is no cancellation
of terms in the expression for $F_t$.

In particular,

$V(F_t)= \min \{ m(g_1) , m(g_2)+1 ,\dots, m(g_r)+r-1 \},$

\noindent   for some $r$ such that $r \leq m(g_1)$ (because each $m(g_k) >
1$), and

$V(F_z)= \min \{ m(f)-1, m(g_1)-1+1 , m(g_2)-1+2 , \dots,  m(g_r)-1+r \}$

$\hskip11mm = \min \{ m - 1, V(F_t) \}.$

\medskip

Why can we affirm this, while Greuel cannot ? It is precisely due to the
separation of the terms in powers of $t$. Any relation of linear
dependence for a fixed $t_0$ between the initial polynomials of the
$g_i$'s will not be preserved as $t$ goes to zero on a linear arc $(uz_0,
ut_0)$. The lowest powers of $t$ will become dominant, leading to our
affirmation in the proof of Proposition 2.2.

\bigskip
{\bf 3. Obtaining equimultiplicity.}

\medskip
{\bf Proposition 3.1.} {\it Let $F(z,t)$ be a $\mu$-constant family of
complex hypersurfaces with isolated singularities at $z = 0$ for each $t$
in a neighbourhood of $0$, of the form $F(z,t) = f(z) + \Sigma_{k=1} t^k
g_k(z)$.  Suppose that the tangent cone of $f$ has an isolated singularity
at $0$. Then the multiplicity  $m(F(.,t))$ is constant as $t$ varies in a
neighbourhood of $0$. }

{\bf Proof.}
Suppose that the tangent cone $\{f_m = 0\}$ (where $f_m$ is the initial
polynomial of $f$) has an isolated singularity. Then $f$ is
semi-homogeneous, in particular semi-quasihomogeneous, and using deep work
of Varchenko [23],  Greuel [8] and O'Shea [16]  independently proved
equimultiplicity.
The special case of homogeneous $f$ was previously treated by Gabrielov
and Koushnirenko [7].
\qed

\medskip

 Motivated by the proof of Proposition 2.1, we could study what happens
in a  family  $F(z,t) = f(z) + tg(z) +t^2h(z)$ with constant Milnor
number  if we take a more general generic curve of the form $ (u^pz_0,
u^qt_0)$ with $z_0 \neq 0, t_0 \neq 0$, and where $p \neq q$ and $p$ and
$q$ are non-negative integers. This turns out not to be fruitful
however.

So we change tactics by choosing an appropriate {\it non-generic} line
segment, whereas the previous results were obtained by choosing suitable
{\it  generic} line segments.

\medskip

{\bf Proposition 3.2.} {\it Let $F(z,t)$ be a $\mu$-constant family of
complex hypersurfaces with isolated singularities at $z = 0$ for each $t$
in a neighbourhood of $0$, of the form $F(z,t) = f(z) + tg(z) + t^2h(z)$. 
Suppose that the singular set of the tangent cone of $\{ f = 0 \}$ is not
contained in the tangent cone of $\{ h = 0 \}$. Then the multiplicity 
$m(F(.,t))$ is constant as $t$ varies in a neighbourhood of $0$.}

\medskip

{\bf Proof.}
By Proposition 2.1, we know that $m(g) \geq m =  m(f)$, and $m(h) \geq m -
1$.

Suppose $m(h) = m - 1.$

For a complex line segment $\gamma(u) = (uz_0,ut_0)$, we have
 $$\hskip-40truemm
 V(\Delta) = V( {{ \vert g+2th \vert} \over {\vert \vert f_z + tg_z +
t^2h_z  \vert \vert }}) $$

 $$\hskip28truemm  = V(g +2th) -  {\text{inf}} \{V(f_{z_i} + tg_{z_i} +
t^2 h_{z_i}) \,  \vert \, i = 1, \dots , n \}. \quad \quad (*)$$

The hypothesis that the singular locus $\Sigma$ of the tangent cone of $\{
f = 0 \}$ is not contained in the tangent cone of $\{ h = 0 \}$ ensures in
particular that $\Sigma$ has complex dimension at least one. Let $z_0$ be
a point of $\Sigma$, distinct from the origin. Then we know that on arcs
$\gamma(u) = (uz_0,ut_0)$, for any choice of non-zero $t_0$,  $v(f_{z_i})
\geq m$ for all $i$.

But since $m(g) \geq m$, and $m(h) = m - 1$, then also $V(tg_{z_i}) \geq
m$ for all $i$, and $V(t^2h_{z_i}) \geq m$, for all $i$. It follows that
$ {\text{inf}} \{V(f_{z_i} + tg_{z_i} + t^2 h_{z_i}) \, \vert \, i = 1,
\dots , n \} \geq m$. Now choose a generic value of $t_0$ for which $V(g
+ 2th) = V(th) = m$.
By $(*)$, $V(\Delta) \leq m - m = 0$, and the L\^e-Saito-Teissier
criterion states that the constancy of the Milnor number implies a
contradiction, i.e.  $m(h) \neq m - 1$. Because $m(h) \geq m - 1$, it
follows that $m(h) \geq m$ and since also $m(g) \geq m$, it follows that
the multiplicity of $F(.,t) = f + tg + t^2h$ is constant as $t$ varies in
a neighbourhood of $0$.
\qed

\bigskip

{\bf Remark 3.3.} Similarly to the argument in the previous proof we can
obtain a contradiction to the hypothesis of constant Milnor number by the
L\^e-Saito-Teissier criterion if $m(g) = m$ and $\Sigma(f_m)$ is not
contained in the tangent cone to $\{g =0 \}$. Thus if a family $F(z,t) =
f(z) + tg(z) + t^2h(z)$ has constant Milnor number, and the singular locus
of the tangent cone of $\{f = 0\}$ is not contained in the tangent cone of
$\{ g = 0 \}$, then $m(g) \geq m+1$.

More generally, the same argument shows that if a family $F(z,t) = f(z) +
tg_1(z) +t^2g_2(z) + t^3g_3(z) + \dots $ has constant Milnor number, and
the singular locus of the tangent cone of $\{f = 0\}$ is not contained in
the tangent cone of $\{ g_k = 0 \}$, then $m(g_k) \geq m-k+2 $.

\bigskip
{\bf Bibliography.}

\medskip

\noindent 1. K. Bekka,  C-r\'egularit\'e et trivialit\'e topologique, {\it
Singularity theory and applications, Warwick 1989 (eds. D. M. Q. Mond and
J. Montaldi)}, Springer Lecture Notes {\bf 1462} (1991), 42-62.

\noindent 2. K. Bekka and D. Trotman, Propri\'et\'e m\'etriques de
familles $\Phi$-radiales de sous-vari\'et\'es diff\'erentiables,  {\it C.
R. Acad. Sci., Paris } {\bf  305} (1987), 389-392.

\noindent 3. K. Bekka and D. Trotman,  Weakly Whitney stratified sets. 
{\it Real and complex singularities
(Proceedings, Sao Carlos 1998, edited by J. W. Bruce and F. Tari)},
{\it Chapman and Hall/CRC},  (2000) 1-15.

\noindent 4. K. Bekka and D. Trotman, Brian\c con-Speder type examples and
the failure of weak Whitney regularity, preprint, 2011.

\noindent 5. J. Brian\c con and J.-P. Speder, La trivialit\'e topologique
n'implique pas les conditions de Whitney, {\it C. R. Acad. Sci., Paris }
{\bf  280} (1975), 365-367.

\noindent 6. J. Brian\c con and J.-P. Speder, Les conditions de Whitney
impliquent $\mu^*$- constant, {\it Annales de l'Institut Fourier,  
Grenoble,} {\bf 26 (2)} (1976), 153-163.

.\noindent 7. A. Gabrielov and A. Kouchnirenko, Description of deformations
with constant Milnor number for homogeneous functions, {\it Funct.
Analysis} {\bf 9} (1975), 67-68.

\noindent 8. G.-M.  Greuel, Constant Milnor number implies constant
multiplicity for quasihomogeneous singularities, {\it Manuscripta Math.}
{\bf 56} (1986), 159-166.

\noindent  9. H. Hironaka, Normal cones of analytic Whitney
stratifications, {\it Publ. Math. IHES.} {\bf 36} (1969), 127-138.

\noindent 10. L\^e D\~ung Tr\'ang and K. Saito, La constance du nombre de
Milnor donne des bonnes stratifications, {\it C. R. Acad. Sci. Paris} {\bf
277} (1973), 793-795.

\noindent 11. L\^e D\~ung Tr\'ang and C. P. Ramanujam, The invariance of
Milnor's number implies the invariance of the topological type, {\it
Amer. J. of Math.} {\bf 98} (1976), 67-78.

\noindent 12. L\^e D\~ung Tr\'ang and B. Teissier, Report on the problem
session, {\it Singularities, Part 2 (Arcata, Calif., 1981),  Proc. Sympos.
Pure Math}., {\bf 40}, Amer. Math. Soc., Providence, R.I. (1983), 105-116.

\noindent 13. J. Milnor, {\it Singular points of complex hypersurfaces},
Annals of Mathematics Studies, Princeton, 1968.

\noindent 14. V. Navarro Aznar, Conditions de Whitney et sections planes,
{\it Inventiones Math. } {\bf  61} (1980), 199-225.

\noindent 15. D. O'Shea, Topologically trivial deformations of isolated
quasihomogeneous hypersurface singularities are equimultiple, {\it Proc.
Amer. Math. Soc.} {\bf 101} (1987), 260-262.

\noindent 16. A. Parusinski, Topological triviality of $\mu$-constant
deformations of type $ f(x) + tg(x)$, {\it Bull. London Math. Soc.} {\bf 
31 } (1999),  686-692.

\noindent 17. B. Teissier, Cycles \'evanescents, sections planes et
conditions de Whitney, {\it Singularit\'es \`a Carg\`ese (ed. F. Pham),
Ast\'erisque} {\bf 7-8} (1973), 285-362.

\noindent 18. B. Teissier, Introduction to equisingularity problems, {\it
Algebraic Geometry (ed. F. Pham), Proc. AMS summer symposia} {\bf 29}
(1975), 593-632 .

 \noindent 19. B. Teissier, Vari\'et\'es polaires II: Multiplicit\'es
polaires, sections planes, et conditions de Whitney, {\it Algebraic
Geometry Proceedings, La Rabida 1981}, Lecture Notes in Math.  {\bf
961}, Springer-Verlag, New York(1982), 314-491.

\noindent 20. D. Trotman, Partial results on the topological invariance
of the multiplicity of a complex hypersurface, {\it S\'eminaire
A'Campo-MacPherson, Universit\'e de Paris 7}, 26 March 1977.

\noindent 21. D. Trotman, {\it \'Equisingularit\'e et conditions de
Whitney}, Th\`ese d'\'Etat, Universit\'e de Paris-Sud, Orsay, 1980.

\noindent 22. D. Trotman and D. Van Straten, Weak Whitney regularity
implies equimultiplicity for families of singular complex analytic
hypersurfaces, preprint, 2011.

\noindent 23 A.N. Varchenko, A lower bound for the codimension of the
stratum $\mu =$  const in terms of the mixed Hodge structure. Vest. Mosk.
Univ. Mat. {\bf 37} (1982), 29-31.

\noindent 24. O. Zariski, Some open questions in the theory of
singularities, {\it Bull. Amer. Math. Soc.} {\bf 77} (1971), 481-491.

\medskip

Laboratoire d'Analyse, Topologie et Probabilit\'es (UMR 7353 du CNRS),

Aix-Marseille Universit\'e,

39 rue Joliot-Curie,
13453 Marseille Cedex 13,
France.

\medskip
plenat\@univ-amu.fr, trotman\@univ-amu.fr

\end